\documentclass[12pt,a4paper,english,reqno]{amsart}
\usepackage[a4paper,footskip=1.5em]{geometry}
\usepackage{amsmath,amssymb,amsthm,mathtools,mathrsfs}
\usepackage{adjustbox,tikz,calc,graphics,babel}
\usepackage{csquotes,enumerate,verbatim}
\usepackage[final]{microtype}
\usepackage[numbers]{natbib}
\usetikzlibrary{shapes.misc,calc,intersections,patterns,decorations.pathreplacing}
\usepackage{hyperref}
\hypersetup{colorlinks=true,linkcolor=blue,citecolor=blue}

\allowdisplaybreaks

\theoremstyle{plain}
\newtheorem*{theorem*}{Theorem}
\newtheorem{theorem}{Theorem}[section]
\newtheorem{lemma}[theorem]{Lemma}
\newtheorem{claim}[theorem]{Claim}

\newtheorem*{claim*}{Claim}

\newtheorem{conjecture}[theorem]{Conjecture}

\theoremstyle{remark}

\def\N{\mathbb{N}}

\def\C{\mathcal}

\let\emptyset\varnothing

\let\originalleft\left
\let\originalright\right
\renewcommand{\left}{\mathopen{}\mathclose\bgroup\originalleft}
\renewcommand{\right}{\aftergroup\egroup\originalright}

\makeatletter
\def\imod#1{\allowbreak\mkern10mu({\operator@font mod}\,\,#1)}
\makeatother

\def\b_sum[#1]{\sum_{i=0}^{r}{\binom{#1}{i}}}

\begin{document}

\title{Approximations to $m$-coloured complete infinite hypergraphs}

\author{Teeradej Kittipassorn}
\address{Department of Mathematical Sciences, University of Memphis, Memphis TN 38152, USA}
\email{t.kittipassorn@memphis.edu}

\author{Bhargav Narayanan}
\address{Department of Pure Mathematics and Mathematical Statistics, University of Cambridge, Wilberforce Road, Cambridge CB3\thinspace0WB, UK}
\email{b.p.narayanan@dpmms.cam.ac.uk}

\date{29 September 2013}
\subjclass[2010]{Primary 05D10; Secondary 05C63, 05C65}

\begin{abstract}
Given an edge colouring of a graph with a set of $m$ colours, we say that the graph is (\emph{exactly}) \emph{$m$-coloured} if each of the colours is used. In 1999, Stacey and Weidl, partially resolving a conjecture of Erickson from 1994, showed that for a fixed natural number $m>2$ and for all sufficiently large $k$, there is a $k$-colouring of the complete graph on $\N$ such that no complete infinite subgraph is exactly $m$-coloured. In the light of this result, we consider the question of how close we can come to finding an exactly $m$-coloured complete infinite subgraph. We show that for a natural number $m$ and any finite colouring of the edges of the complete graph on $\N$ with $m$ or more colours, there is an exactly ${\hat m}$-coloured complete infinite subgraph for some ${\hat m}$ satisfying $|m-{\hat m}|\le \sqrt{m/2} + 1/2$; this is best-possible up to the additive constant. We also obtain analogous results for this problem in the setting of $r$-uniform hypergraphs. Along the way, we also prove a recent conjecture of the second author and investigate generalisations of this conjecture to $r$-uniform hypergraphs.
\end{abstract}

\maketitle

\section{Introduction}
The classical problem of Ramsey theory is to find a large monochromatic structure in a larger coloured structure; for a host of results, see~\citep{Graham1990}. On the other hand, the objects of interest in anti-Ramsey theory are large `rainbow coloured' or `totally multicoloured' structures; see, for example, the paper of Erd{\H{o}}s,  Simonovits and S{\'o}s~\cite{Erdos1975}. Between these two ends of the spectrum, one could consider the question of finding structures which are coloured with exactly $m$ different colours: this was first done by Erickson~\citep{Erickson1994} and this is the line of enquiry that we pursue here.

Our notation is standard. Thus, following Erd{\H{o}}s, for a set $X$, we write  $X^{(r)}$ for the family of all subsets of $X$ of cardinality $r$; equivalently, $X^{(r)}$ is the complete $r$-uniform hypergraph on the vertex set $X$. We write $[n]$ for $\{1,\dots,n\}$, the set of the first $n$ natural numbers.  We denote a surjective map $f$ from a set $X$ to another set $Y$ by $f: X \twoheadrightarrow Y$. By a \emph{colouring} of a hypergraph, we mean a colouring of the edges of the hypergraph unless we specify otherwise.

Let $\Delta:\N^{(r)}\twoheadrightarrow[k]$ be a surjective $k$-colouring of the edges of the complete $r$-uniform hypergraph on the natural numbers. We say that a subset $X\subset\N$ is (\emph{exactly}) $m$\emph{-coloured} if $\Delta( X^{(r)} )$, the set of values attained by $\Delta$ on the edges induced by $X$, has size exactly $m$. Let $\gamma_{\Delta}(X)$, or $\gamma(X)$ in short, denote the size of the set $\Delta( X^{(r)})$; in other words, every set $X$ is $\gamma(X)$-coloured. In this paper, we shall study for fixed $r$ and large $k$, the set of values $m$ for which there exists an infinite $m$-coloured set with respect to a $k$-colouring $\Delta:\N^{(r)}\twoheadrightarrow[k]$. Let us mention as an aside that it is also interesting to study what happens when we wish to find finite $m$-coloured sets, or allow colourings which use infinitely many colours; we refer the reader to~\citep{BK2013} for results of this flavour. With our goal of finding infinite $m$-coloured sets in mind, let us define, for a $k$-colouring $\Delta:\N^{(r)}\twoheadrightarrow[k]$, the set
\[
\C{F}_{\Delta}=\left\{ \gamma_\Delta(X)  :  X\subset\N \mbox{ such that } X \mbox{ is infinite}\right\} .
\]
Clearly, $k\in\C{F}_{\Delta}$ as $\Delta$ is surjective, and Ramsey's Theorem tells us that $1\in\C{F}_{\Delta}$. Erickson~\citep{Erickson1994} noted that a fairly straightforward application of Ramsey's Theorem enables one to show that $2\in\C{F}_{\Delta}$ for any $k$-colouring $\Delta$ of $\N^{(r)}$ with $k\ge 2$.  He also conjectured that with the exception of $1,2\mbox{ and }k$, no other elements are guaranteed to be in $\C{F}_{\Delta}$ (even in the case of graphs) and that if $k>m>2$, then there is a $k$-colouring $\Delta$ of $\N^{(2)}$ such that $m\notin\C{F}_{\Delta}$. Stacey and Weidl~\citep{Stacey1999}, partially resolving this conjecture, showed using a probabilistic construction that there is a constant $C_{m}$ such that if $k>C_{m}$, then there is a $k$-colouring $\Delta$ of $\N^{(2)}$ such that $m\notin\C{F}_{\Delta}$.

Since an exactly $m$-coloured complete infinite subhypergraph is not guaranteed to exist, we are naturally led to the question of whether we can find a complete infinite subhypergraph that is exactly ${\hat m}$-coloured for some ${\hat m}$ close to $m$. In this paper, we establish the following result.

\begin{theorem}
\label{approx13-approxthm}
Fix a positive integer $r \ge 2$. For any $k$-colouring $\Delta:\N^{(r)}\twoheadrightarrow[k]$ and any natural number $ m \le k$, there exists an ${\hat m} \in \C{F}_\Delta $ such that
\[|m-{\hat m}| \le c_r m^{1-1/r} + O(m^{1-2/r}),\]
where $c_r = r/(2(r!)^{1/r})$.
\end{theorem}

Theorem~\ref{approx13-approxthm} is tight up to the $O(m^{1-2/r})$ term. To see this, let $k=\binom{n}{r}+1$ for some $n \in \N$. We consider the `small-rainbow colouring' $\Delta$ which colours all the edges induced by $[n]$ with $\binom{n}{r}$ distinct colours and all the remaining edges with the one colour that has not been used so far. In this case, we see that $\C{F}_\Delta  = \{ \binom{i}{r}+1 : i \le n\}$. Now let $m=(\binom{l}{r} + \binom{l+1}{r} + 2)/2$ for some natural number $l$ such that $l < n$. It is not difficult to check that $|m-{\hat m}| \ge \binom{l}{r-1}/2$ for each ${\hat m} \in \C{F}_\Delta $; also, it is clear that $\binom{l}{r-1}/2 = (c_r-o(1)) m^{1-1/r}$.

In the case of graphs where $r=2$, Theorem~\ref{approx13-approxthm} tells us that for any finite colouring of the edges of the complete graph on $\N$ with $m$ or more colours, there is an exactly ${\hat m}$-coloured complete infinite subgraph for some ${\hat m}$ satisfying $|m-{\hat m}|\le \sqrt{m/2} + O(1)$; a careful analysis of the proof of Theorem~\ref{approx13-approxthm} in this case allows us to replace  the $O(1)$ term with an explicit constant, $1/2$.

We know from Theorem~\ref{approx13-approxthm} that $\C{F}_\Delta $ cannot contain very large gaps. Another natural question
we are led to ask is if there are any sets, and in particular, intervals that $\C{F}_\Delta$ is guaranteed to intersect. Making this more precise, the second author conjectured, see~\citep{Bhargav2013}, that the small-rainbow colouring described above is extremal for graphs in the following sense.

\begin{conjecture}
\label{approx13-choose2}
Let $\Delta:\N^{(2)}\twoheadrightarrow[k]$ be a $k$-colouring of the complete graph on $\N$ and suppose $n$ is a natural number such that $k>\binom{n}{2}+1$. Then $ \C{F}_\Delta \cap (\binom{n}{2}+1,\binom{n+1}{2} +1 ] \not= \emptyset$.
\end{conjecture}

In this paper, we shall prove this conjecture. There are two natural generalisations of this conjecture to $r$-uniform hypergraphs which are equivalent to Conjecture~\ref{approx13-choose2} in the case of graphs.

The first comes from considering small-rainbow colourings; indeed we can ask whether $\C{F}_\Delta  \cap I_{r,n} \neq \emptyset$ when $k > \binom{n}{r}+1$, where $I_{r,n}$ is the interval $( \binom{n}{r}+1, \binom{n+1}{r}+1]$.

The second comes from considering a different family of colourings which we call `small-set colourings'. Let $k = \b_sum[n]$ and consider the surjective $k$-colouring $\Delta$ of $\N^{(r)}$ defined by $\Delta(e) = e \cap [n]$. Note that in this case, $\C{F}_\Delta  = \lbrace \b_sum[j]: j\le n \rbrace$. Consequently, we can ask whether $\C{F}_\Delta  \cap J_{r,n} \neq \emptyset$ when $k > \b_sum[n-1]$, where $J_{r,n}$ is the interval $( \b_sum[n-1], \b_sum[n]]$.

Note that both these questions are identical when $r = 2$. Indeed, $\binom{n}{2} +\binom{n}{1} + \binom{n}{0} = \binom{n+1}{2}+1$, so $I_{2,n} = J_{2,n}$.

We shall demonstrate that the correct generalisation is the former. We shall first prove that the answer to the first question is in the affirmative, provided $n$ is sufficiently large.

\begin{theorem}
\label{approx13-interval-I}
For every $r\ge2$, there exists a natural number $n_r \ge r-1$ such that for any natural number $n\ge n_r$ and any $k$-colouring $\Delta:\N^{(r)}\twoheadrightarrow[k]$  with $k>\binom{n}{r}+1$, $\C{F}_\Delta \cap I_{r,n} \not= \emptyset$.
\end{theorem}

Using a result of Baranyai~\citep{Baranyai1975} on factorisations of uniform hypergraphs, we shall exhibit an infinite family of colourings that answer the second question negatively for every $r\ge 3$.

\begin{theorem}
\label{approx13-interval-J}
For every $r\ge3$, there exist infinitely many values of $n$ for which there exists a $k$-colouring $\Delta:\N^{(r)}\twoheadrightarrow[k]$ with $k > \b_sum[n-1]$ such that $\C{F}_\Delta \cap J_{r,n} = \emptyset$.
\end{theorem}

The rest of this paper is organised as follows. In the next section, we shall prove Theorems~\ref{approx13-approxthm},~\ref{approx13-interval-I} and~\ref{approx13-interval-J} and deduce Conjecture~\ref{approx13-choose2} from the proof of Theorem~\ref{approx13-interval-I}. We then conclude by mentioning some open problems.

\section{Proofs of the main results}

We start with the following lemma which we shall later use to prove both Theorems~\ref{approx13-approxthm} and~\ref{approx13-interval-I}.

\begin{lemma}
\label{approx13-mainlem}
Let $m\ge2$ be an element of $\C{F}_\Delta $. Then there exists a natural number $a=a(m,\Delta)$ such that
\begin{enumerate}
\item $\b_sum[a] \ge m$, and
\item $\C{F}_\Delta \cap [m - \min(\sum_{i=0}^{r-1}\binom{a-1}{i},{r(m-1)/a}), m) \not= \emptyset$.
\end{enumerate}
Futhermore, if
\[m=\sum_{i=t+1}^{r}\binom{a}{i}+s+1\]
for some $s\ge 0$ and $0\le t+1\le r$, then
\[\C{F}_\Delta \cap \left[\sum_{i=t+1}^{r}\binom{a-1}{i}+\left(1-\frac{t}{a}\right)s+1,m\right) \not= \emptyset.\]
\end{lemma}

\begin{proof} We start by establishing the following claim.

\begin{claim}
\label{approx13-X-A}
There is an infinite $m$-coloured set $X\subset \N$ with a finite subset $A\subset X$ such that
\begin{enumerate}
	\item \label{approx13-one}the colour of every edge of $X$ is determined by its intersection with $A$, i.e., if $e_1\cap A=e_2\cap A$, then $\Delta(e_1)=\Delta(e_2)$, and
	\item \label{approx13-two}$\gamma(X\setminus\{v\})<m$ for all $v\in A$.
\end{enumerate}
\end{claim}
\begin{proof}
To see this, let $W\subset\N$ be an infinite $m$-coloured set. For each colour $c\in\Delta(W^{(r)})$, pick an edge $e_c$ in $W$ of colour $c$ and let $A = \bigcup_c e_c$ be the set of vertices incident to these edges. So $A\subset W$ is a finite $m$-coloured set. Let $A_1,A_2,\dots,A_l$ be an enumeration of the subsets of $A$ of size at most $r$. Note that this is the complete list of possible intersections of an edge with $A$. We now define a descending sequence of infinite sets $B_0\supset B_1\supset\dots\supset B_l$ as follows. Let $B_0=W\setminus A$. Having defined the infinite set $B_{i-1}$, we induce a colouring of the  $(r-|A_i|)$-tuples $T$ of $B_{i-1}$, by giving $T$ the colour of the edge $A_i\cup T$. By Ramsey's Theorem, there is an infinite monochromatic subset $B_i\subset B_{i-1}$ with respect to this induced colouring, so the edges of $A\cup B_i$ whose intersection with $A$ is $A_i$ have the same colour.

Hence, $X=A\cup B_l$ is an infinite $m$-coloured set satisfying property (\ref{approx13-one}). Now, if we have a vertex $v\in A$ such that $\gamma(X\setminus\{v\})=m$, we delete $v$ from $A$. We repeat this until we are left with an $m$-coloured set $X$ satisfying (\ref{approx13-one}) and (\ref{approx13-two}).
\end{proof}

Let $X$ and $A$ be as guaranteed by Claim~\ref{approx13-X-A}. Note that $A$ is nonempty since $m\ge2$. We shall prove the lemma with $a(m,\Delta) = |A|$. From the structure of $X$ and $A$, we note that $\b_sum[a] \ge m$. That
\[\C{F}_\Delta \cap \left[m - \min\left(\sum_{i=0}^{r-1}\binom{a-1}{i},\frac{r(m-1)}{a}\right), m\right) \not= \emptyset\] is a consequence of the following claim.

\begin{claim}
\label{approx13-X1-X2}
There exist infinite sets $X_1, X_2 \subset X$ such that $m-\sum_{i=0}^{r-1}\binom{a-1}{i}\le\gamma(X_1)<m$ and $m-r(m-1)/a\le\gamma(X_2)<m$.
\end{claim}

\begin{proof}
Let $X_1=X\setminus\{v\}$ for any $v\in A$. We know from Claim~\ref{approx13-X-A} that $\gamma(X_1)<m$. We shall now prove that $\gamma(X_1)\ge m-\sum_{i=0}^{r-1}\binom{a-1}{i}$; that is, the number of colours lost by removing $v$ from $X$ is at most $\sum_{i=0}^{r-1}\binom{a-1}{i}$. Since the colour of an edge is determined by its intersection with $A$, the number of colours lost is at most the numbers of subsets of $A$ containing $v$ of size at most $r$, which is precisely $\sum_{i=0}^{r-1}\binom{a-1}{i}$.

Next, we shall prove that there is a subset $X_2\subset X$ such that $m-r(m-1)/a\le\gamma(X_2)<m$. Let $A=\{v_1,v_2,\dots,v_a\}$ and let
	\[
	C_i=\Delta\left(X^{(r)}\right)\setminus \Delta\left((X\setminus\{v_i\})^{(r)}\right)
	\]
be the set of colours lost by removing $v_i$ from $X$; since $\gamma(X\setminus\{v_i\})<m$ for all $v_i\in A$, it follows that $C_i\not=\emptyset$. For each colour $c \in \Delta ( X^{(r)} )$, pick an edge $e_c$ of colour $c$, and let $A_c= e_c\cap A$; in particular, we take $A_{c_\emptyset}=\emptyset$, where $c_\emptyset$ is the colour corresponding to an empty intersection with $A$. Since every edge of colour $c\in C_i$ contains $v_i$, we double count the number of times a colour is counted in the sum $\sum_{i=1}^a{|C_i|}$ to obtain
	\[
	\sum_{i=1}^a{|C_i|}\le \sum_{c\not=c_\emptyset}{|A_c|}\le r(m-1),
	\]
so there exists an $i$ such that $0<|C_i|\le{r(m-1)/a}$; the claim follows by taking $X_2=X\setminus\{v_i\}$.
\end{proof}

We finish the proof of the lemma by establishing the following claim.
\begin{claim}
\label{approx13-rep}
If we can write $m=\sum_{i=t+1}^{r}\binom{a}{i}+s+1$ for some $s\ge 0$ and $0\le t+1\le r$, then
\[\C{F}_\Delta \cap \left[\sum_{i=t+1}^{r}\binom{a-1}{i}+\left(1-\frac{t}{a}\right)s+1,m\right) \not= \emptyset.\]
\end{claim}
\begin{proof}
As in the proof of Claim~\ref{approx13-X1-X2}, for each colour $c \in \Delta ( X^{(r)} )$, pick an edge $e_c$ of colour $c$, and let $A_c= e_c\cap A$; in particular, let $A_{c_\emptyset}=\emptyset$. We know from Claim~\ref{approx13-X-A} that edges of $X$ of distinct colours cannot have the same intersection with $A$. Consequently, all the $A_c$ are distinct subsets of $A$, each of size at most $r$. Hence,
	\[\sum_{c\not=c_\emptyset}{|A_c|}\le\sum_{i=t+1}^{r} i \binom{a}{i} +ts.\]
Arguing as in the proof of Claim~\ref{approx13-X1-X2}, we conclude that there exists a vertex $v\in A$ such that the number of colours lost by removing $v$ from $X$ is at most $(\sum_{i=t+1}^{r} i \binom{a}{i} +ts)/a$. Therefore,
	\begin{align*}
		\gamma(X\setminus \{v\})
		&\ge m-\frac{1}{a}\left(\sum_{i=t+1}^{r} i \binom{a}{i} +ts\right)\\
		&=m-\left(\sum_{i=t+1}^{r} \binom{a-1}{i-1} +\frac{ts}{a}\right)\\
		&=\sum_{i=t+1}^{r}  \binom{a-1}{i}+\left(1-\frac{t}{a}\right)s+1,
	\end{align*}
so it follows that
\[ \C{F}_\Delta \cap \left[\sum_{i=t+1}^{r}\binom{a-1}{i}+\left(1-\frac{t}{a}\right)s+1,m\right) \not= \emptyset. \qedhere \]
\end{proof}

The lemma now follows from Claims~\ref{approx13-X-A},~\ref{approx13-X1-X2} and~\ref{approx13-rep}.
\end{proof}

Having established Lemma~\ref{approx13-mainlem}, it is easy to deduce both Theorem~\ref{approx13-approxthm} and~\ref{approx13-interval-I} from the lemma.

\begin{proof}[Proof of Theorem~\ref{approx13-approxthm}]
Let $t = m+c_r m^{1-1/r}$. We may assume that $m > r^r/r!$ since otherwise $m= O(1)$ and there is nothing to prove. Also, if $t \ge k$, then the result follows easily by taking $\hat m = k$ so we may assume that $t < k$. Let ${\hat t}$ be the smallest element of $\C{F}_\Delta $ greater than $t$. Applying Lemma~\ref{approx13-mainlem} to ${\hat t}$, we find an ${\hat m} \in \C{F}_\Delta $ such that ${\hat m} \le t$ and
\[{\hat m} \ge {\hat t} -  \min\left(\sum_{i=0}^{r-1}\binom{a-1}{i},\frac{r({\hat t}-1)}{a}\right)\]
for some natural number $a$.  Now if $a \ge (r!m)^{1/r} > r$,  then
\[{\hat m} \ge {\hat t} - \frac{r({\hat t}-1)}{a} \ge {\hat t}\left(1 - \frac{r}{a}\right)\ge t\left(1 - \frac{r}{a}\right), \]
so it follows that ${\hat m} \ge m - c_r m^{1-1/r} - O(m^{1-2/r})$. If $a < (r!m)^{1/r}$ on the other hand, then using the fact that
\begin{align*}
{\hat m} &\ge {\hat t} - \sum_{i=0}^{r-1}\binom{a-1}{i} \\
&\ge t - \frac{a^{r-1}}{(r-1)!} - O(a^{r-2}) \\
&\ge t - \frac{(r!m)^{1-1/r}}{(r-1)!} - O(m^{1-2/r}),
\end{align*} it follows once again that ${\hat m} \ge m - c_r m^{1-1/r} - O(m^{1-2/r})$.
\end{proof}

\begin{proof}[Proof of Theorem~\ref{approx13-interval-I}]
If $k \le \binom{n+1}{r}+1$, we are done since $k\in\C{F}_\Delta $. So suppose that $k >\binom{n+1}{r} +1$. Let $m$ be the smallest element of $\C{F}_\Delta $ such that $m>\binom{n+1}{r}+1$; hence, $\C{F}_\Delta \cap (\binom{n+1}{r}+1,m)=\emptyset$. Now, since $m\ge2$, there exists by Lemma~\ref{approx13-mainlem}, a natural number $a$ such that
\[\C{F}_\Delta \cap \left[m - \frac{r(m-1)}{a}, \binom{n+1}{r}+1\right] \not= \emptyset.\]
To prove the theorem, it is sufficient to show that $m - {r(m-1)/a}>\binom{n}{r}+1$. We know from Lemma~\ref{approx13-mainlem} that $\b_sum[a] \ge m > \binom{n+1}{r}+1$. If $n$ is sufficiently large, we must have $a\ge n$.

If $a\ge n+1$, then
\begin{align*}
m-\frac{r(m-1)}{a} &= (m-1)\left(1-\frac{r}{a}\right)+1 \\
&> \binom{n+1}{r}\left(1-\frac{r}{n+1}\right)+1\\
&=\binom{n}{r}+1
\end{align*}
since $m> \binom{n+1}{r}+1$ and $n\ge r-1$.

We now deal with the case $a=n$. First, we write $m=\binom{n}{r}+\binom{n}{r-1}+s+1$. Since $m> \binom{n+1}{r}+1$ and $\binom{n}{r}+\binom{n}{r-1} = \binom{n+1}{r}$, we see that $s>0$. By Lemma~\ref{approx13-mainlem}, it follows that
\[\C{F}_\Delta \cap \left[\binom{n}{r}+\left(1-\frac{r-2}{n}\right)s+1, m\right) \not= \emptyset.\]
Since $n\ge r-1$ and $s>0$, the result follows.
\end{proof}

A careful inspection of the proof of Theorem~\ref{approx13-interval-I} shows that when $r=2$, the statement holds for all $n\in \N$. We hence obtain a proof of Conjecture~\ref{approx13-choose2}. By constructing a sequence of highly structured subgraphs, the second author~\citep{Bhargav2013} proved that for any $k$-colouring $\Delta:\N^{(2)}\twoheadrightarrow[k]$ with $k\ge\binom{n}{2}+1$ for some natural number $n$, $| \C{F}_\Delta  | \ge n$; our proof of Conjecture~\ref{approx13-choose2} gives a short proof of this lower bound. Theorem~\ref{approx13-interval-I} also yields a generalisation of this lower bound for $r$-uniform hypergraphs, albeit with a constant additive error term (which depends on $r$).

We now turn to the proof of Theorem~\ref{approx13-interval-J}. We will need a result of Baranyai's~\citep{Baranyai1975} which states that the set of edges of the complete $r$-uniform hypergraph on $l$ vertices can be partitioned into perfect matchings when $r \,\vert\,l$.

\begin{proof}[Proof of Theorem~\ref{approx13-interval-J}]
We shall show that if $n$ is sufficiently large and $(r-1)\,\vert\, (n+1)$, then there is a surjective $k$-colouring $\Delta$ of $\N^{(r)}$ with $k>\b_sum[n-1]$ and $\C{F}_\Delta \cap J_{r,n} = \emptyset$. We shall define a colouring of $\N^{(r)}$ such that the colour of an edge $e$ is determined by its intersection with a set $A$ of size $n+1$, say $A = [n+1]$. Let $\C{B}$ be the family of all subsets of $A$ of size at most $r$. For $B\in \C{B}$, we denote the colour assigned to all the edges $e$ such that $e\cap A = B$ by $c_B$.

To define our colouring, we shall construct a partition $\C{B}=\C{B}_1\cup\C{B}_2$ with $\emptyset\in\C{B}_2$. Then for every $B\in\C{B}_2$, we set $c_B$ to be equal to $c_\emptyset$. Finally, we take the colours $c_B$ for $B\in\C{B}_1$ to all be distinct and different from $c_\emptyset$. Hence, the number of colours used is $k=|\C{B}_1|+1$. It remains to construct this partition of $\C{B}$.

Since $(r-1)\,\vert\, (n+1)$, by Baranyai's theorem there exists an ordering
\[B_1,B_2,\dots,B_{\binom{n+1}{r-1}}\] of the subsets of $A$ of size $r-1$ such that for all $0\le t \le \binom{n}{r-2}$, the family
\[\left\lbrace B_{\left(\frac{n+1}{r-1}\right)t+1},B_{\left(\frac{n+1}{r-1}\right)t+2},\dots,B_{\left(\frac{n+1}{r-1}\right)(t+1)} \right\rbrace\]
is a perfect matching. Let $\C{B}_1=\{B_1,B_2,\dots,B_s\}\cup\{B\in\C{B}:|B|=r\}$, where
\[s=\b_sum[n]-\binom{n+1}{r};\]
our colouring is well defined because $0\le s\le \binom{n+1}{r-1}$ for all sufficiently large $n$. Observe that
\[k=|\C{B}_1|+1=\binom{n+1}{r}+s+1=\b_sum[n]+1.\]
We shall show that the second largest element of $\C{F}_\Delta $ is at most $\b_sum[n-1]$. Note that any $X\subset \N$ with $\gamma(X) < k$ cannot contain $A$. As before, let $C_i$ be the set of colours lost by removing $i\in A$ from $\N$, i.e.,
\[C_i = \Delta\left(\N^{(r)}\right) \setminus\Delta\left((\N\setminus\{i\})^{(r)}\right).\]
We shall complete the proof by showing that $k-|C_i|\le \b_sum[n-1]$ for all $i \in A$.

Note that our construction ensures that $||C_i|-|C_j||\le 1$ for all $i,j \in A$. Now, observe that
	\[\sum_{i=1}^{n+1}{|C_i|}= \sum_{B\in \C{B}_1}{|B|}= r\binom{n+1}{r}+(r-1)s,\]
so $|C_i|\ge (r\binom{n+1}{r}+(r-1)s)/(n+1)-1$ for all $i\in A$. It is then easily verified using Pascal's identity that when $r \ge 4$ and $n$ is sufficiently large,
 	\begin{align*}
		k-|C_i|&\le \left(\binom{n+1}{r}+s+1\right)-\frac{1}{n+1}\left(r\binom{n+1}{r}+(r-1)s\right)+1 \\
		&= \binom{n}{r}+\left(1-\frac{r-1}{n+1}\right)s+2\\
		&= \binom{n}{r}+\left(1-\frac{r-1}{n+1}\right)\left(\b_sum[n]-\binom{n+1}{r}\right)+2\\
		&\le \b_sum[n-1];
	\end{align*}
the last inequality above is deduced by comparing the coefficients of the polynomials in the inequality.

When $r = 3$, it is easy to check that $s = n+1$, so $s$ is divisible by $(n+1)/(r-1) = (n+1)/2$. Consequently, in this case, $|C_i| = |C_j|$ for $i,j \in A$. Hence,
	\begin{align*}
		k-|C_i|&\le \left(\binom{n+1}{3}+s+1\right)-\frac{1}{n+1}\left(3\binom{n+1}{3}+2s\right) \\
		&= \binom{n}{3}+\left(1-\frac{2}{n+1}\right)(n+1)+1\\
		&= \sum_{i=0}^{3}{\binom{n-1}{i}}.
	\end{align*}
This completes the proof.
\end{proof}

\section{Conclusion}

We conclude by mentioning two open problems. We proved that for any $k$-colouring $\Delta:\N^{(r)}\twoheadrightarrow[k]$ and every sufficiently large natural number $n$, $\C{F}_\Delta \cap I_{r,n} \not= \emptyset$ provided $k>\binom{n}{r}+1$. A careful analysis of our proof shows that the result holds when $n \ge (5/2 + o(1))r$; we chose not to give details to keep the presentation simple. However, we suspect that the result should hold as long as $n\ge r-1$ but a proof eludes us.

To state the next problem, let us define
\[
\psi_r(k)=\min_{\Delta:\N^{(r)}\twoheadrightarrow[k]}|\C{F}_{\Delta}|.
\]
A consequence of Theorem~\ref{approx13-interval-I} is that $\psi_r (k) \ge (r!k)^{1/r} - O(1)$. Turning to the question of upper bounds for $\psi_r$, the small-rainbow colouring shows that the lower bound that we get from Theorem~\ref{approx13-interval-I} is tight infinitely often, i.e., when $k$ is of the form $\binom{n}{r} + 1$ for some $n \in \N$. When $k$ is not of this form, there are two obvious ways of generalising the small-rainbow colouring: we could replace the rainbow coloured clique in our construction either with a disjoint union of cliques or with a clique along with a pendant vertex attached to some subset of the vertices of the clique. However, both these obvious generalisations of the small-rainbow colouring fail to give us good upper bounds for $\psi_r (k)$ for a general $k \in \N$. The second author proved~\citep{Bhargav2013} using rainbow colourings of complete bipartite graphs that
\[ \psi_2 (k) =O \left(\frac{k}{(\log{\log{k}})^{\delta}(\log{\log{\log{k}}})^{3/2}}\right)\]
for almost all natural numbers $k$ and some absolute constant $\delta >0$. The same construction can be extended to show that $\psi_r(k) = o(k)$ for almost all natural numbers $k$. It would be very interesting to decide if, in fact, $\psi_r(k)=o(k)$ for all $k \in \N$.

\section*{Acknowledgements} Some of the research in this paper was carried out while the authors were visitors at Microsoft Research, Redmond. We are grateful to Yuval Peres and the other members of the Theory Group at Microsoft Research for their hospitality.

We would also like to thank Tomas Ju\v{s}kevi\v{c}ius for helpful discussions.

\bibliographystyle{amsplain}
\bibliography{approx_col}

\end{document}